\documentclass[12pt]{amsart}

\usepackage{amsmath}
\usepackage{amsfonts}
\usepackage{amssymb}
\usepackage{amsthm}
\usepackage{mathrsfs}
\usepackage{hyperref}
\usepackage{enumerate}
\usepackage{cleveref}
\usepackage{mathrsfs}
\usepackage[top=0.9in, bottom=1.15in, left=1.1in, right=1.1in]{geometry}

\newtheorem{theorem}{Theorem}[section]

\newtheorem{proposition}[theorem]{Proposition}
\newtheorem{corollary}[theorem]{Corollary}
\newtheorem{definition}[theorem]{Definition}
\newtheorem*{remark}{Remark}

\Crefname{conjecture}{Conjecture}{Conjectures}

\theoremstyle{plain}

\theoremstyle{plain}

\renewcommand{\P}{\mathbb{P}}
\newcommand{\Z}{\mathbb{Z}}
\newcommand{\Q}{\mathbb{Q}}

\newcommand{\C}{\mathbb{C}}

\numberwithin{equation}{section}

\author{Robert Schneider}
\address{Department of Mathematics and Computer Science\newline
Emory University\newline
400 Dowman Dr., W401\newline
Atlanta, GA 30322}
\email{robert.schneider@emory.edu}
\title{Arithmetic of partitions and the $q$-bracket operator}

\begin{document}
\begin{abstract}
We present a natural multiplicative theory of integer partitions (which are usually considered in terms of addition), and find many theorems of classical number theory arise as particular cases of extremely general combinatorial structure laws. We then see that the relatively recently-defined $q$-bracket operator $\left<f\right>_q$, studied by Bloch--Okounkov, Zagier, and others for its quasimodular properties, plays a deep role in the theory of partitions, quite apart from questions of modularity. Moreover, we give an explicit formula for the coefficients of  $\left<f\right>_q$ for any function $f$ defined on partitions, and, conversely, give a partition-theoretic function whose $q$-bracket is a given power series.  
\end{abstract}

\bibliographystyle{amsplain}
\maketitle

\section{Introduction: the \texorpdfstring{$q$}{Lg}-bracket operator}\label{section1}

In a groundbreaking paper of 2000 \cite{BlochOkounkov}, Bloch and Okounkov introduced the $q$-bracket operator $\left<f\right>_q$ of a function $f$ defined on the set of integer partitions, and showed that the $q$-bracket can be used to produce quasimodular forms. A recent paper \cite{Zagier} by Zagier examines the $q$-bracket operator from various enlightening perspectives, and finds broader classes of quasimodular forms arising from its application. This study is inspired by Zagier's treatment, as well as by ideas of Alladi--Erd\H{o}s \cite{AlladiErdos}, Andrews \cite{Andrews}, and Fine \cite{Fine}.

We fix a few notations and concepts in order to proceed. Let $\mathcal P$ denote the set of integer partitions $\lambda = (\lambda_1,\lambda_2,\dots, \lambda_r)$, $\lambda_1 \geq \lambda_2 \geq \cdots \geq \lambda_r \geq 1$, including the ``empty partition" $\emptyset$. We call the number of parts of $\lambda$ the \textit{length} $\ell(\lambda) := r$ of the partition. We call the number being partitioned the \textit{size} $|\lambda|:= \lambda_1 + \lambda_2 + \cdots + \lambda_r$ of $\lambda$. We assume the conventions $\ell(\emptyset) := 0$ and $|\emptyset|:= 0$. Let us write $\lambda \vdash n$ to indicate $\lambda$ is a partition of $n$ (i.e., $|\lambda| = n$), and we will allow a slight abuse of notation to let $\lambda_i \in \lambda$ indicate $\lambda_i$ is one of the parts of $\lambda$.

Furthermore, taking $z\in\C$ and $q := e^{2\pi i\tau}$ with $\tau$ in the upper half-plane $\mathbb H$ (thus $|q|<1$), we define the $q$-Pochhammer symbol $(z;q)_{*}$ by $(z;q)_0:=1, \  (z;q)_n := \prod_{k=0}^{n-1}(1-zq^k)$ for $n\geq 1$, and $(z;q)_{\infty} := \lim_{n\to \infty}(z;q)_n$. We assume absolute convergence everywhere. 


Now we may give the definition of the $q$-bracket operator of Bloch and Okounkov. 



\begin{definition}\label{qbracket} 
We define the {\it $q$-bracket} $\left<f\right>_q$ of a function $f:\mathcal P \to \C$ to be the quotient
\[
\left<f\right>_q:=\frac{\sum_{\lambda \in \mathcal P}f(\lambda)q^{|\lambda|}}{\sum_{\lambda \in \mathcal P}q^{|\lambda|}}\in \C[[q]]
.\]
We take the resulting power series to be indexed by partitions, unless otherwise specified.

\end{definition}

\begin{remark} 
Definition \ref{qbracket} extends the range of the $q$-bracket somewhat; the operator is defined in \cite{BlochOkounkov} and \cite{Zagier} to be a power series in $\Q[[q]]$ instead of $\C[[q]]$, as those works take $f:\mathcal P\to\Q$. 
We may write the $q$-bracket in equivalent forms that will prove useful here:  
\begin{equation}
\left<f\right>_q=\left(q;q\right)_{\infty}\sum_{\lambda \in \mathcal P}f(\lambda)q^{|\lambda|}=\left(q;q\right)_{\infty}\sum_{n=0}^{\infty}q^n\sum_{\lambda \vdash n}f(\lambda)
\end{equation}

\end{remark}

While computationally the operator boils down to multiplying a power series by $(q;q)_{\infty}$, conceptually the $q$-bracket represents a sort of weighted average of the function $f$ over all partitions. Zagier gives an interpretation of the $q$-bracket as the ``expectation value of an observable $f$ in a statistical system whose
states are labelled by partitions'' \cite{Zagier}. Such sums over partitions are ubiquitous in statistical mechanics, quantum physics, and string theory \cite{Okounkov, Zwiebach}. We will keep in the backs of our minds the poetic feeling that the partition-theoretic structures we encounter are, somehow, part of the fabric of physical reality. 

We begin our study by considering the $q$-bracket of a prominent statistic in partition theory, the {\it rank} function $\operatorname{rk}(\lambda)$ of Freeman Dyson \cite{Dyson}, defined by
$$\operatorname{rk}(\lambda):=\operatorname{lg}(\lambda)-\ell(\lambda)$$
where we let $\operatorname{lg}(\lambda)$ denote the {\it largest part} of the partition (similarly, we write $\operatorname{sm}(\lambda)$ for the {\it smallest part}). Noting that $\sum_{\lambda\vdash n}\operatorname{rk}(\lambda)=1$ if $n=0$ (i.e., if $\lambda=\emptyset$) and is equal to $0$ otherwise, then 
$$\sum_{\lambda\in\mathcal P}\operatorname{rk}(\lambda)q^{|\lambda|}=\sum_{n=0}^{\infty}q^n\sum_{\lambda\vdash n}\operatorname{rk}(\lambda)=1.$$
Therefore we have that
\begin{equation}\label{rank}
\left<\operatorname{rk}\right>_q=(q;q)_{\infty}.
\end{equation}
We see by comparison with the Dedekind eta function $\eta(\tau):=q^{\frac{1}{24}}(q;q)_{\infty}$ that $\left<\operatorname{rk}\right>_q$ is very nearly a weight-$1/2$ modular form. 

Now, recall the weight-$2k$ Eisenstein series central to the theory of modular forms \cite{Ono}
\begin{equation}\label{Eisenstein}
E_{2k}(\tau)=1-\frac{4k}{B_{2k}}\sum_{n=1}^{\infty}\sigma_{2k-1}(n)q^n,
\end{equation}
where $k\geq 1$, $B_j$ denotes the $j$th Bernoulli number, and $\sigma_{*}$ is the classical sum-of-divisors function. It is not hard to see the $q$-bracket of the ``size'' function
\[
\left<|\cdot|\right>_q=-q\frac{\frac{d}{dq}{(q;q)_{\infty}}}{(q;q)_{\infty}}=\frac{1-E_2(\tau)}{24}
\]
is essentially quasimodular; the series $E_2(\tau)$ is the prototype of a quasimodular form.  

The near-modularity of the $q$-bracket of basic partition-theoretic functions is among the operator's most fascinating features. Bloch--Okounkov give a recipe for constructing quasimodular forms using $q$-brackets of shifted symmetric polynomials \cite{BlochOkounkov}. Zagier expands on their work to find infinite families of quasimodular $q$-brackets, including families that lie outside Bloch and Okounkov's methods \cite{Zagier}. Griffin--Jameson--Trebat-Leder build on these methods to find $p$-adic modular and quasimodular forms as well \cite{GJTL}. While it appears at first glance to be little more than convenient shorthand, the $q$-bracket notation identifies---induces, even---intriguing classes of partition-theoretic  phenomena. 

In this study, we give an exact formula for the coefficients of $\left<f\right>_q$ for any function $f:\mathcal P\to\C$. We also answer the converse problem, viz. for an arbitrary power series $\hat{f}(q)$ we give a function $F$ defined on $\mathcal P$ such that $\left<F\right>_q=\hat{f}(q)$ exactly. The main theorems appear in Section \ref{3}. 

Along the way, we establish a simple, general multiplicative theory of integer partitions, which specializes to many fundamental results in classical number theory. Moreover, we adopt a certain philosophical position here:  Multiplicative number theory in $\Z$ is a special case of vastly general combinatorial laws, one out of an infinity of parallel number theories in a partition-theoretic multiverse. 

It turns out the $q$-bracket operator plays a surprisingly natural role in this multiverse.

\  
\  
\  
\

\section{Multiplicative arithmetic of partitions}\label{Multiplicative theory of partitions}

We introduce one more statistic alongside the length $\ell(\lambda)$ and size $|\lambda|$ of a partition, that we call the \textit{integer} of $\lambda$. 

\begin{definition}
We define the ``integer'' of a partition $\lambda$, written $n_{\lambda}$, to be the product of the parts of $\lambda$, i.e., we take 
$$n_{\lambda}:=\lambda_1\lambda_2\cdots \lambda_r.$$ 
We adopt the convention $n_{\emptyset}:=1$ (it is an empty product).
\end{definition}

It may not seem to be a very natural statistic---after all, partitions are defined additively with no straightforward connection to multiplication---but this so-called ``integer'' of a partition shows up in partition-theoretic formulas scattered throughout the literature \cite{Andrews, Fine, Robert}, and will prove to be important to the theory indicated here as well. 

Pushing further in the multiplicative direction, we define a simple, intuitive multiplication operation on the elements of $\mathcal P$.

\begin{definition}\label{productdef}
We define the \textit{product} $\lambda \lambda'$ of two partitions $\lambda,\lambda'\in\mathcal P$ as the multi-set union of their parts listed in weakly decreasing order, e.g. $(5,3,2,2)(4,2,1,1)=(5,4,3,2,2,2,1,1)$. The empty partition $\emptyset$ serves as the multiplicative identity. 
\end{definition}


Then it makes sense to write $\lambda^2:=\lambda \lambda, \lambda^3:=\lambda \lambda \lambda$, and so on. It is easy to see that we have the following relations:
\begin{align*}\label{log relations}
&n_{\lambda\lambda'}=n_{\lambda}n_{\lambda'}, &&n_{\lambda^a}=n_{\lambda}^a\\
&\ell(\lambda\lambda')=\ell(\lambda)+\ell(\lambda'), &&\ell(\lambda^a)=a\cdot\ell(\lambda)\\
&|\lambda\lambda'|=|\lambda|+|\lambda'|, &&|\lambda^a|=a |\lambda|
\end{align*}
Note that length and size both resemble logarithms, which O. Beckwith pointed out to the author \cite{Olivia}.

We also define division in $\mathcal P$.
\begin{definition}\label{divisiondef}
We will say a partition $\delta$ \textit{divides} $\lambda$ and will write $\delta | \lambda$, if all of the parts of $\delta$ are also parts of $\lambda$, including multiplicity, e.g. we write $(4,2,1,1)|(5,4,3,2,2,2,1,1)$. When $\delta | \lambda$ we might also discuss the quotient $\lambda / \delta\in\mathcal P$ (or $\frac{\lambda}{\delta}$) formed by deleting the parts of $\delta$ from $\lambda$.
\end{definition}

\begin{remark}
That is, $\delta$ is a ``sub-partition'' of $\lambda$. Note that the empty partition $\emptyset$ is a divisor of every partition $\lambda$.
\end{remark}
 
So it makes sense to write $\lambda^0:=\lambda/\lambda=\emptyset$. We also have relations such as these:
$$n_{\lambda / \lambda'}=\frac{n_{\lambda}}{n_{\lambda'}},\  \  \ell(\lambda / \lambda')=\ell(\lambda)-\ell(\lambda'),\  \  |\lambda / \lambda'|=|\lambda|-|\lambda'|
$$



On analogy to the prime numbers in classical arithmetic, the partitions into one part (e.g. $(1),(3),(4)$) are both prime and irreducible under this simple multiplication. The analog of the Fundamental Theorem of Arithmetic is trivial: of course, every partition may be uniquely decomposed into its parts. Thus we might rewrite a partition $\lambda$ in terms of its ``prime'' factorization
$\lambda=(a_1)^{m_1} (a_2)^{m_2}...(a_t)^{m_t}$, where $a_1 > a_2 > ... > a_t \geq 1$ are the distinct numbers appearing in $\lambda$ such that $a_1=\operatorname{lg}(\lambda)$ (the largest part of $\lambda$), $a_t=\operatorname{sm}(\lambda)$ (the smallest part), and $m_i$ denotes the multiplicity of $a_i$ as a part of $\lambda$. Clearly, then, we have 
\begin{equation}
n_{\lambda}=a_1^{m_1} a_2^{m_2}\cdots a_t^{m_t}.
\end{equation} 
We note in passing that we also have a dual formula for the ``integer'' $n_{\lambda^*}$ of the conjugate $\lambda^*$ of $\lambda$, written in terms of $\lambda$, viz.
\begin{equation}
n_{\lambda^*}=M_1^{a_1-a_2}M_2^{a_2-a_3}\cdots M_{t-1}^{a_{t-1}-a_t}M_t^{a_t},
\end{equation}
where $M_k:=\sum_{i=1}^{k}m_i$ (thus $M_t=\ell(\lambda)$), which is evident from the Young diagrams of $\lambda$ and $\lambda^*$.

Fundamental classical concepts such as coprimality, greatest common divisor, least common multiple, etc., apply with exactly the same meanings in the partition-theoretic setting, if one replaces ``prime factors of a number'' with ``parts of a partition'' in the classical definitions. 

\begin{remark}
If $\mathcal P' \subseteq \mathcal P$ is an infinite subset of $\mathcal P$ closed under partition multiplication and division, then the multiplicative theory presented in this study still holds when the relations are restricted to $\mathcal P'$ (this condition is sufficient, at least). In particular, a partition ideal of order 1 in the sense of Andrews \cite{Andrews}---such as the set of partitions into prime parts---is closed under partition multiplication and division, and conforms to the structures noted here.
\end{remark}

\  
\  
\  
\

\section{Partition-theoretic analogs of classical functions}\label{muphisigma}

A number of important functions from classical number theory have partition-theoretic analogs, giving rise to nice summation identities that generalize their classical counterparts. We see examples of this phenomenon in the author's study of partition zeta functions \cite{Robert}.

One of the most fundamental classical arithmetic functions, related to factorization of integers, is the M{\"o}bius function. It turns out there is a natural partition-theoretic analog  of $\mu$.

\begin{definition}\label{moebiusdef}
For $\lambda\in\mathcal P$ we define a partition-theoretic M{\"o}bius function $\mu(\lambda)$ as follows:
$$
\mu(\lambda):= \left\{
        \begin{array}{ll}
            0 & \text{if $\lambda$ has any part repeated,}\\
            (-1)^{\ell(\lambda)} & \text{otherwise}
        \end{array}
    \right.
$$
\end{definition}

If $\lambda$ is a partition into prime parts, the above definition gives the classical M{\"o}bius function $\mu(n_{\lambda})$. Just as in the classical case, we have by inclusion-exclusion the following, familiar relation.

\begin{proposition}
Summing $\mu(\delta)$ over the divisors $\delta$ of $\lambda\in\mathcal P$, we have
$$
\sum_{\delta|\lambda}\mu(\delta)= \left\{
        \begin{array}{ll}
            1 & \text{if $\lambda=\emptyset$,}\\
            0 & \text{otherwise.}
        \end{array}
    \right.
$$
\end{proposition}

Furthermore, we have a partition-theoretic generalization of the M{\"o}bius inversion formula, which is proved along the lines of proofs of the classical formula.
\begin{proposition}\label{mobinv}
For a function $f:\mathcal P \to \C$ we have the equivalence
$$
F(\lambda)=\sum_{\delta|\lambda}f(\delta)\  \Longleftrightarrow\  f(\lambda)=\sum_{\delta|\lambda}F(\delta)\mu(\lambda / \delta)
.$$

\end{proposition}

\begin{remark}
K. Alladi has considered similar partition M{\"o}bius function relations in lectures and unpublished work \cite{Alladi2015}. 
\end{remark}

In classical number theory, M{\"o}bius inversion is often used in conjunction with order-of-summation swapping principles for double summations. These have an obvious partition-theoretic generalization as well, reflected in the following identity.

\begin{proposition}\label{sumswap}
For functions $f,g:\mathcal P \to \C$ we have the formula
$$\sum_{\lambda\in\mathcal P} f(\lambda) \sum_{\delta | \lambda} g(\delta) = \sum_{\lambda\in\mathcal P} g(\lambda) \sum_{\gamma\in\mathcal P} f(\lambda \gamma).
$$
\end{proposition}


The preceding propositions will prove useful in the next section, to understand the coefficients of the $q$-bracket operator. 

The classical M{\"o}bius function has a close companion in the Euler phi function $\varphi(n)$, also known as the totient function, which counts the number of natural numbers less than $n$ that are coprime to $n$. This sort of statistic does not seem meaningful in the partition-theoretic frame of reference, as there is not generally a well-defined greater- or less-than ordering of partitions. However, if we sidestep this business of ordering and counting for the time being, we find it is possible to define a partition analog of $\varphi$ that is naturally compatible with the identities above, as well as with classical identities involving the Euler phi function. 

Recall that $n_{\lambda}$ denotes the ``integer'' of $\lambda$, i.e., the product of its parts.

\begin{definition}\label{phidef}
For $\lambda\in\mathcal P$ we define a partition-theoretic phi function $\varphi(\lambda)$ by
$$
\varphi(\lambda):=\  n_{\lambda}\prod_{\substack{\lambda_i\in\lambda\\ \text{without}\\  \text{repetition}}}(1-\lambda_i^{-1})
,$$
where the product is taken over only the distinct numbers composing $\lambda$, that is, the parts of $\lambda$ without repetition.
\end{definition}

If $\lambda$ is a partition into prime parts, this reduces to the classical $\varphi(n_{\lambda})$. Clearly if $1\in\lambda$ then $\varphi(\lambda)=0$, which is a bit startling by comparison with the classical phi function that never vanishes. This phi function filters out partitions containing 1's. (This is a convenient property: partitions containing 1's can create convergence issues in partition zeta functions \cite{Robert} and other partition-indexed series.)

As with the M{\"o}bius function above, the partition-theoretic $\varphi(\lambda)$ yields generalizations of many classical expressions. For instance, there is a familiar-looking divisor sum, which is proved along classical lines.

\begin{proposition}\label{phisum}
We have that
$$
\sum_{\delta|\lambda}\varphi(\delta)=n_{\lambda}
.$$
\end{proposition}  

We also find a partition analog of the well-known relation connecting the phi function and the M{\"o}bius function, by directly expanding the right-hand side of Definition \ref{phidef}.

\begin{proposition}\label{phimoeb}
We have the identity
$$
\varphi(\lambda)=n_{\lambda}\sum_{\delta|\lambda}\frac{\mu(\delta)}{n_{\delta}}
.$$
\end{proposition}  

Combining the above relations, we arrive at a nicely balanced identity. 

\begin{proposition}\label{phimoeb2}
For $f:\mathcal P\to\C$ let $\widetilde{f}(\lambda):=\sum_{\gamma\in \mathcal P}f(\lambda \gamma)$. Then we have
$$
\sum_{\lambda\in\mathcal P}\frac{f(\lambda)\varphi(\lambda)}{n_{\lambda}}
=\sum_{\lambda\in\mathcal P}\frac{\widetilde{f}(\lambda)\mu(\lambda)}{n_{\lambda}}
.$$
\end{proposition}

\begin{remark}
Replacing $\mathcal P$ with the set $\mathcal P_{\P}$ of partitions into prime parts (the so-called ``prime partitions''), then for $n=n_{\lambda}, k=n_{\gamma}$, the sum $\widetilde{f}$ takes the form $\widetilde{f}(n)=\sum_{k\geq 1}f(nk)$ and we obtain the following identity, which must be known classically:
$$
\sum_{n=1}^{\infty}\frac{f(n)\varphi(n)}{n}
=\sum_{n=1}^{\infty}\frac{\widetilde{f}(n)\mu(n)}{n}
$$
\end{remark}


A number of other important arithmetic functions have partition-theoretic analogs, too, such as the sum-of-divisors function $\sigma_a$.

\begin{definition}\label{sigma}
For $\lambda\in\mathcal P, a\in\Z_{\geq 0}$, we define the function 
$$
\sigma_a(\lambda):=\sum_{\delta|\lambda}n_{\delta}^a,
$$
with the convention $\sigma(\lambda):=\sigma_1(\lambda)$.
\end{definition}
If $\lambda\in\mathcal P_{\P}$ this becomes the classical function $\sigma_a(n_{\lambda})$. One might wonder about ``perfect partitions'' or other analogous phenomena related to $\sigma_a$ classically. This partition sum-of-divisors function will come into play in the next section.

%
%
%

We reiterate, these familiar-looking identities not only mimic classical theorems, they fully generalize the classical cases. The propositions above all specialize to their classical counterparts when we restrict our attention to the set $\mathcal P_{\P}$ of prime partitions; then, as a rule-of-thumb, we just replace partitions with their ``integers'' in the formulas (other parameters may need to be adjusted appropriately). This is due to the bijective correspondence between natural numbers and $\mathcal P_{\P}$ noted in \cite{Robert}, following Alladi and Erd\H{o}s \cite{AlladiErdos}: the set of ``integers'' of prime partitions (including $n_{\emptyset}$) is precisely the set of positive integers $\Z^+$, by the Fundamental Theorem of Arithmetic. Yet prime partitions form a narrow slice, so to speak, of the set $\mathcal P$ over which these general relations hold sway. 

Many well-known laws of classical number theory arise as special cases of underlying partition-theoretic structures.

\  
\  
\  
\

\section{Role of the \texorpdfstring{$q$}{Lg}-bracket} \label{3}

We return now to the $q$-bracket operator of Bloch--Okounkov, which we recall from Definition \ref{qbracket}. The $q$-bracket arises naturally in the multiplicative theory outlined above. To see this, take $F(\lambda):=\sum_{\delta|\lambda}f(\delta)$ for $f:\mathcal P \to\C$. It follows from Proposition \ref{sumswap} that
\begin{align*}
\sum_{\lambda\in\mathcal P} F(\lambda)q^{|\lambda|} & =\sum_{\lambda\in\mathcal P} q^{|\lambda|} \sum_{\delta | \lambda} f(\delta) = \sum_{\lambda\in\mathcal P} f(\lambda) \sum_{\gamma\in\mathcal P} q^{|\lambda \gamma|}
\\
&=\sum_{\lambda\in\mathcal P} f(\lambda) \sum_{\gamma\in\mathcal P} q^{|\lambda|+|\gamma|}
=\left(\sum_{\lambda\in\mathcal P} f(\lambda)q^{|\lambda|}\right) \left(\sum_{\gamma\in\mathcal P} q^{|\gamma|}\right)
.
\end{align*}


Observing that the rightmost sum above is equal to $(q;q)_{\infty}^{-1}$, then by comparison with Definition \ref{qbracket} of the $q$-bracket operator, we arrive at the two central theorems of this study. Together they give a type of $q$-bracket inversion, converting partition divisor sums into partition power series, and vice versa. 

\begin{theorem}\label{thm1}
For an arbitrary function $f:\mathcal P \to\C$, if
$$
F(\lambda)=\sum_{\delta|\lambda}f(\delta)\
$$
then
$$
\left<F\right>_q=\sum_{\lambda\in\mathcal P}f(\lambda)q^{|\lambda|}.
$$
\end{theorem}
 
In the converse direction, we can also write down a simple function whose $q$-bracket is a given power series indexed by partitions.

\begin{theorem}\label{thm1.5}
Consider an arbitrary power series  of the form
$$\sum_{\lambda\in\mathcal P}f(\lambda)q^{|\lambda|}.
$$
Then we have the function $F:\mathcal P\to\C$ given by
$$F(\lambda)=\sum_{\delta|\lambda}f(\delta),
$$  
such that $\left<F\right>_q=\sum_{\lambda\in\mathcal P} f(\lambda)q^{|\lambda|}.$
\end{theorem}



\begin{remark}
These theorems still hold if everywhere (including in the $q$-bracket definition) we replace $\mathcal P$ with a subset $\mathcal P'\subseteq\mathcal P$ closed under partition multiplication and division.
\end{remark}

We wish to apply Theorems \ref{thm1} and \ref{thm1.5} to examine the $q$-brackets of partition-theoretic analogs of classical functions introduced in Section \ref{muphisigma}. 

Recall Definition \ref{sigma} of the sum of divisors function $\sigma_a(\lambda)$. Then $\sigma_0(\lambda)=\sum_{\delta | \lambda}1$ counts the number of partition divisors (i.e., sub-partitions) of $\lambda \in\mathcal P$, much as in the classical case. It is immediate from Theorem \ref{thm1} that
\begin{equation}\label{sigmabracket}
\left<\sigma_0\right>_q=(q;q)_{\infty}^{-1}.
\end{equation}
If we note that $(q;q)_{\infty}$ is also a factor of the $q$-bracket on the left-hand side, we can see as well
\begin{equation}\label{sigmabracket2}
\sum_{\lambda\in\mathcal P}\sigma_0(\lambda)q^{|\lambda|}=(q;q)_{\infty}^{-2}
.
\end{equation}
Remembering also from Equation \ref{rank} the identity $\left<\operatorname{rk}\right>_q=(q;q)_{\infty}$, so far we have seen a few instances of interesting power series connected to powers of $(q;q)_{\infty}$ via the $q$-bracket operator. We return to this point below. 


Now let us recall the handful of partition divisor sum identities from Section \ref{muphisigma} involving the partition-theoretic functions $\varphi,\sigma_a$, and the ``integer of a partition'' function $n_{*}$. Theorem \ref{thm1} reveals that these three functions form a close-knit family, related through (double) application of the $q$-bracket.

\begin{corollary}\label{qbracketsystem}
We have the pair of identities 
\begin{align*}\label{system}
\left<\sigma\right>_q &=\sum_{\lambda \in \mathcal P}n_{\lambda}q^{|\lambda|}
\\
\left<n_{*}\right>_q &=\sum_{\lambda \in \mathcal P}\varphi(\lambda)q^{|\lambda|}
.\end{align*}
\end{corollary}



The coefficients of $\left<\sigma\right>_q$ are of the form $n_{*}$; applying the $q$-bracket a second time to the function $n_{*}$ gives us the rightmost summation, whose coefficients are the values of $\varphi$. In fact, it is evident that this operation of applying the $q$-bracket more than once might be continued indefinitely; thus we feel the need to introduce a new notation, on analogy to differentiation.

\begin{definition}\label{diffndef}
If we apply the $q$-bracket repeatedly, say $n\geq 0$ times, to the function $f$, we denote this operator by $\left<f\right>_q^{(n)}$. We define $\left<f\right>_q^{(n)}$ by the equation
$$\left<f\right>_q^{(n)}:=(q;q)_{\infty}^n\sum_{\lambda \in \mathcal P}f(\lambda)q^{|\lambda|}\in\C[[q]].$$
\end{definition}

\begin{remark}
It follows from the definition above that $\left<f\right>_q^{(0)}=\sum_{\lambda \in \mathcal P}f(\lambda)q^{|\lambda|},\  \left<f\right>_q^{(1)}=\left<f\right>_q$. 
\end{remark}

Theorem \ref{thm1.5} gives us a converse construction as well, allowing us to write down a partition divisor sum whose $q$-bracket is a given power series. Then we define an inverse ``antibracket'', analogous to the antiderivative in calculus.

\begin{definition}

We call $F:\mathcal P\to \C$ a ``$q$-{\it antibracket}'' of $f$ if $\left<F\right>_q= \sum_{\lambda \in \mathcal P}f(\lambda) q^{|\lambda|}$.
\end{definition}

 
The act of finding an antibracket might be carried out repeatedly as well. 
We define a canonical class of $q$-antibrackets related to $f$ by extending Definition \ref{diffndef} to allow for negative values of $n$.

\begin{definition}\label{antidef}
If we repeatedly divide the power series $\sum_{\lambda\in\mathcal P} f(\lambda)q^{|\lambda|}$ by $(q;q)_{\infty}$, say $n > 0$ times, we notate this operator as 
$$\left<f\right>_q^{(-n)}:=(q;q)_{\infty}^{-n}\sum_{\lambda \in \mathcal P}f(\lambda)q^{|\lambda|}\in\C[[q]].$$
We take the resulting power series to be indexed by partitions, unless otherwise specified. We call the function on $\mathcal P$ defined by the coefficients of $\left<f\right>_q^{(-1)}$ the ``canonical $q$-antibracket'' of $f$ (or sometimes just ``the antibracket'').
\end{definition}

Taken together, Definitions \ref{diffndef} and \ref{antidef} describe an infinite family of $q$-brackets and antibrackets. The following identities give an example of such a family (and of the use of these new bracket notations).

\begin{corollary}\label{compact}

Corollary \ref{qbracketsystem} can be written more compactly as 
$$
\left<\sigma\right>_q^{(2)}=\left<n_{*}\right>_q^{(1)}=\left<\varphi\right>_q^{(0)}
.$$ 
We can also condense Corollary \ref{qbracketsystem} by writing
$$
\left<\sigma\right>_q^{(0)}=\left<n_{*}\right>_q^{(-1)}=\left<\varphi\right>_q^{(-2)}
.$$
\end{corollary}

Both of the compact forms above preserve the essential message of Corollary \ref{qbracketsystem}, that these three partition-theoretic functions are directly connected through the $q$-bracket operator, or more concretely (and perhaps more astonishingly), simply through multiplication or division by powers of $(q;q)_{\infty}$. 

Along similar lines, we can encode Equations \ref{rank}, \ref{sigmabracket}, and \ref{sigmabracket2} in a single statement, noting an infinite family of power series that contains $\left<\operatorname{rk}\right>_q$ and $\left<\sigma_0\right>_q$.

\begin{corollary}\label{ranksigma}
For $n\in\Z$, we have the family of $q$-brackets
\[
\left<\operatorname{rk}\right>_q^{(n)}
=\left<\sigma_0\right>_q^{(n+2)}\\
=(q;q)_{\infty}^n
.\]
\end{corollary}

\begin{remark}
Here we see the $q$-bracket connecting with modularity properties. For instance, another member of this family is $\left<\operatorname{rk}\right>_q^{(24)}=q^{-1}\Delta (\tau)$, where $\Delta$ is the important modular discriminant function having Ramanujan's tau function as its coefficients \cite{Ono}.
\end{remark}

The identities above worked out easily because we knew in advance what the coefficients of the $q$-brackets should be, due to the divisor sum identities from Section \ref{muphisigma}. Theorems \ref{thm1} and \ref{thm1.5} provide a recipe for turning partition divisor sums $F$ into coefficients $f$ of power series, and vice versa. 

However, generally a function $F:\mathcal P\to\C$ is not given as a sum over partition divisors. If we wished to write it in this form, what function $f:\mathcal P\to\C$ would make up the summands? In classical number theory this question is answered by the M{\"o}bius inversion formula; indeed, we have the partition-theoretic analog of this formula in Proposition \ref{mobinv}. 

Recall the ``divided by'' notation $\lambda / \delta$ from Definition \ref{divisiondef}. Then we may write the function $f$ (and thus the coefficients of $\left<F\right>_q$) explicitly using partition M{\"o}bius inversion. 

\begin{theorem}\label{thm2}
The $q$-bracket of the function $F:\mathcal P\to\C$ is given explicitly by
$$
\left<F\right>_q=\sum_{\lambda \in \mathcal P}f(\lambda)q^{|\lambda|}
,
$$
where the coefficients can be written in terms of $F$ itself:
$$
f(\lambda)=\sum_{\delta|\lambda}F(\delta)\mu(\lambda / \delta)
$$
\end{theorem}


We already know from Theorem \ref{thm1.5} that the coefficients of the canonical antibracket of $f$ are written as divisor sums over values of $f$. Thus, much like $\operatorname{rk}(\lambda)$ in Corollary \ref{ranksigma}, every function $f$ defined on partitions can be viewed as the generator, so to speak, of the (possibly infinite) family of power series $\left<f\right>_q^{(n)}$ for $n\in\Z$, whose coefficients can be written in terms of $f$ as $n$-tuple sums of the shape $\sum_{\delta_1|\lambda}\sum_{\delta_2|\delta_1}...\sum_{\delta_n|\delta_{n-1}}$ constructed by repeated application of the above theorems.
%
%
%
%
%

This suggests the following useful fact. 

\begin{corollary}\label{two-way}
If two power series are members of the family $\left<f\right>_q^{(n)}\  \left(n\in\Z\right),$ then the coefficients of each series can be written explicitly in terms of the coefficients of the other. 
\end{corollary}

\  
\  
\  
\

\section{The \texorpdfstring{$q$}{Lg}-antibracket and coefficients of power series over \texorpdfstring{$\Z_{\geq 0}$}{Lg}}\label{Applications of the $q$-antibracket}

Theorems \ref{thm1}, \ref{thm1.5}, and \ref{thm2} together provide a two-way map between the coefficients of families of power series indexed by partitions. In this section, we address the question of computing the antibracket (loosely speaking) of coefficients indexed not by partitions, but by natural numbers $\Z_{\geq 0}$ as usual. 

We remark immediately that coefficients of the form $c_n$ may be expressed in terms of partitions in a number of ways, which are generally not equivalent in this framework of $q$-brackets and antibrackets. For instance, as M. Jameson pointed out to the author \cite{Marie}, we can always write 
\begin{equation}\label{MarieIdentity}
c_n=\sum_{\lambda\vdash n} \frac{c_{\scriptscriptstyle{|\lambda|}}}{p(|\lambda|)}.
\end{equation}
Thus there is more than one function $F:\mathcal P\to\C$ with $\left<F\right>_q=\sum_{n=0}^{\infty} c_n q^n$ for a given sequence $c_n$ of coefficients. (Note that if $F_1, F_2$ are two such functions, we still have $\sum_{\lambda \vdash n}F_1(\lambda)=\sum_{\lambda \vdash n}F_2(\lambda)$, by the equality of coefficients indexed by integers.)

Here we construct an antibracket $F$ using ideas developed in the preceding sections. There are three classes of power series of the form $\sum_{n=0}^{\infty} c_n q^n$ that we examine: (1) the coefficients $c_n$ are given as sums $\sum_{\lambda \vdash n}$ over partitions of $n$; (2) the coefficients $c_n$ are sums $\sum_{d | n}$ over divisors of $n$; and (3) the coefficients $c_n$ are an arbitrary sequence of complex numbers. 

The class (1) above is already given by Theorem \ref{thm1}; to keep this section relatively self-contained, we rephrase the result here. 

\begin{corollary}\label{class1}
For $c_n=\sum_{\lambda \vdash n}f(\lambda)$ we can write
\begin{equation*}
\sum_{n=0}^{\infty}c_n q^n=\sum_{\lambda\in\mathcal P}f(\lambda)q^{|\lambda|}.
\end{equation*}
Then we have a function $F(\lambda)=\sum_{\delta|\lambda}f(\delta)$ such that 
$
\left<F\right>_q=\sum_{n=0}^{\infty}c_n q^n
$.
\end{corollary}

Thus the power series of class (1) are already in a form subject to the $q$-bracket machinery detailed in the previous section. The class (2) with coefficients of the form $\sum_{d|n}$ is a little more subtle. We introduce a special subset $\mathcal P_{=}$ which bridges sums over partitions and sums over the divisors of natural numbers. 

\begin{definition}\label{=def}
We define the subset $\mathcal P_{=}\subseteq\mathcal P$ to be the set of partitions into equal parts, that is, whose parts are all the same positive number, e.g. $(1),(1,1),(4,4,4)$. We make the assumption $\emptyset\notin\mathcal P_{=}$, as the empty partition has no positive parts.
\end{definition}

The divisors of $n$ correspond exactly (in two different ways) to the set of partitions of $n$ into equal parts, i.e., partitions of $n$ in $\mathcal P_{=}$. For example, compare the divisors of $6$
$$
1,2,3,6
$$
with the partitions of $6$ into equal parts 
$$
(6),(3,3),(2,2,2),(1,1,1,1,1,1)
.$$ 
Note that for each of the above partitions $(a,a,...,a)\vdash 6$, we have that $a\cdot\ell\left((a,a,...,a)\right)$ $=6$. We see from this example that for any $n\in\Z^+$ we can uniquely associate each divisor $d|n$ to a partition $\lambda\vdash n,\lambda\in\mathcal P_{=},$ by taking $d$ to be the length of $\lambda$. (Alternatively, we could identify the divisor $d$ with  $\operatorname{lg}(\lambda)$ or $\operatorname{sm}(\lambda)$, as defined above, which of course are the same in this case. We choose here to associate divisors to $\ell(\lambda)$ as length is a universal characteristic of partitions, regardless of the structure of the parts.) 

By the above considerations, it is clear that
\begin{equation}\label{=prop}
\sum_{d|n}f(d)=\sum_{\substack{\lambda\vdash n\\  \lambda\in\mathcal P_{=}}}f\left(\ell(\lambda)\right).
\end{equation}
This leads us to a formula for the coefficients of a power series of the class (2) discussed above.  

\begin{corollary}\label{class2}
For $c_n=\sum_{d|n}f(d)$ we can write
\begin{equation*}\label{4.1}
\sum_{n=0}^{\infty}c_n q^n=\sum_{\lambda\in\mathcal P_{=}}f\left(\ell(\lambda)\right)q^{|\lambda|}.
\end{equation*}
Then we have a function 
$$F(\lambda)=\sum_{\substack{\delta|\lambda \\  \delta\in\mathcal P_{=}}}f\left(\ell(\delta)\right)$$ such that $
\left<F\right>_q=\sum_{n=0}^{\infty}c_n q^n
$.
\end{corollary}

The completely general class (3) of power series with arbitrary coefficients $c_n\in\C$ follows right away from Corollary \ref{class2} by classical M\"{o}bius inversion.

\begin{corollary}\label{class3}
For $c_n\in\C$ we can write
\begin{equation*}
\sum_{n=0}^{\infty}c_n q^n=\sum_{\lambda\in\mathcal P_{=}}q^{|\lambda|}\sum_{d|\ell(\lambda)}c_d\  \mu\left(\frac{\ell(\lambda)}{d}\right).
\end{equation*}
Then we have a function 
$$F(\lambda)=\sum_{\substack{\delta|\lambda \\  \delta\in\mathcal P_{=}}}\sum_{d|\ell(\delta)}c_d\  \mu\left(\frac{\ell(\delta)}{d}\right)$$ such that $
\left<F\right>_q=\sum_{n=0}^{\infty}c_n q^n
$.
\end{corollary}

\begin{remark} 
We point out an alternative expression for sums of the shape of $F(\lambda)$ in this corollary, that can be useful for computation. If we write out the factorization of a partition $\lambda=$ $(a_1)^{m_1}(a_2)^{m_2}...$ $(a_t)^{m_t}$ as in Section \ref{Multiplicative theory of partitions}, a divisor of $\lambda$ lying in $\mathcal P_{=}$ must be of the form $(a_i)^m$ for some $1\leq i \leq t$ and $1 \leq m \leq m_i$. Then for any function $\phi$ defined on $\Z^+$ we see 
\begin{equation}
\sum_{\substack{\delta|\lambda \\  \delta\in\mathcal P_{=}}}\phi\left(\ell(\delta)\right)=\sum_{i=1}^{t}\sum_{j=1}^{m_i}\phi\left(\ell((a_i)^j)\right)=\sum_{i=1}^{t}\sum_{j=1}^{m_i}\phi(j).
\end{equation}
We could also use Equation \ref{MarieIdentity} together with Corollary \ref{class1} to address series of the class (3); however, computation of the partition function $p(*)$ presents a practical challenge.
\end{remark}

Given the ideas developed above, we can now pass between $q$-brackets and arbitrary power series, summed over either natural numbers or partitions.

\  
\  
\  
\

\section{Applications of the \texorpdfstring{$q$}{Lg}-bracket and \texorpdfstring{$q$}{Lg}-antibracket}\label{applications}
We close this report by briefly illustrating some of the methods of the previous sections through two examples.

\subsection{Sum of divisors function}
In classical number theory, for $a\geq 0$ the divisor sum $\sigma_a(n):=\sum_{d|n}d^a$ is particularly important to the theory of modular forms; as seen in Equation \ref{Eisenstein}, for odd values of $a$, power series  of the form
$$
\sum_{n=0}^{\infty}\sigma_a(n)q^n
$$
comprise the Fourier expansions of Eisenstein series \cite{Ono}, which are the building blocks of modular and quasimodular forms. As a straightforward application of Corollary \ref{class2} following directly from the definition of $\sigma_a(n)$, we give a function $\mathcal S_a$ defined on partitions whose $q$-bracket is the power series above.

\begin{corollary}\label{sigmasum}
We have the partition-theoretic function 
$$\mathcal S_a(\lambda):=\sum_{\substack{\delta|\lambda \\  \delta\in\mathcal P_{=}}}\ell(\delta)^a
$$
such that 
$$\left<\mathcal S_a\right>_q=\sum_{n=0}^{\infty}\sigma_a(n) q^n
.$$
\end{corollary}

\begin{remark}
Zagier gives a different function $S_{2k-1}(\lambda)=\sum_{\lambda_i\in\lambda}\lambda_i^{2k-1}$ (the moment function) that also has the $q$-bracket $\sum \sigma_{2k-1}(n) q^n$ \cite{Zagier}. This is an example of the non-uniqueness of antibrackets of functions defined on natural numbers, as noted previously.
\end{remark}

Thus we see the $q$-bracket operator brushing up against modularity, once again.


\subsection{Reciprocal of the Jacobi triple product}

We turn our attention now to another fundamental object in the subject of modular forms. Let $j(z;q)$ denote the classical \textit{Jacobi triple product} \cite{Berndt}
\begin{equation}\label{jtp}
j(z;q):=(z;q)_{\infty}(z^{-1}q;q)_{\infty}(q;q)_{\infty}.
\end{equation}

The reciprocal of the triple product 
$$j(z;q)^{-1}=\sum_{\lambda\in\mathcal P}j_z(\lambda)q^{|\lambda|}$$
is interesting in its own right. For instance, $j(z;q)^{-1}$ plays a role not unlike the role played by $(q;q)_{\infty}$ in the $q$-bracket operator, for the Appell--Lerch sum $m(x,q,z)$ important to the study of mock modular forms \cite{HickersonMortenson}.

Our goal will be to derive a formula for the coefficients $j_z(\lambda)$ above. If we multiply $j(z;q)^{-1}$ by $(1-z)$ to cancel the pole at $z=1$, it behaves nicely under the action of the $q$-bracket. Let us write
\begin{equation}\label{normalized}
(1-z)j(z;q)^{-1}=\frac{1}{(zq;q)_{\infty}(z^{-1}q;q)_{\infty}(q;q)_{\infty}}=\sum_{\lambda\in\mathcal P}J_z(\lambda)q^{|\lambda|}.
\end{equation}

Let $\operatorname{crk}(\lambda)$ denote the {\it crank} of a partition, an important partition-theoretic statistic conjectured by Dyson \cite{Dyson} and discovered almost half a century later by Andrews and Garvan \cite{AndrewsGarvan} in work related to Ramanujan congruences. Crank is not unlike Dyson's rank, but is a bit more complicated (thus we do not give it explicitly). 

We define $M(n,m)$ to be the number of partitions of $n$ having crank equal to $m\in\Z$; then the Andrews--Garvan \textit{crank generating function} $C(z;q)$ is given by
\begin{equation}\label{crankgen}
C(z;q):=\frac{(q;q)_{\infty}}{(zq;q)_{\infty}(z^{-1}q;q)_{\infty}}=\sum_{n=0}^{\infty}M_z(n)q^n,
\end{equation}
where we set
\begin{equation}\label{M}
M_z(n):=\sum_{\lambda\vdash n}z^{\operatorname{crk}(\lambda)}=\sum_{m=-\infty}^{\infty}M(n,m)z^m.
\end{equation}

The function $C(z;q)$ has deep connections. When $z=1$, Equation \ref{crankgen} reduces to Euler's partition generating function formula \cite{Berndt}. For $\zeta\neq 1$ a root of unity, $C(\zeta;q)$ is a modular form, and Folsom--Ono--Rhoades show the crank generating function to be related to the theory of quantum modular forms \cite{FOR}.  

Comparing Equations \ref{normalized} and \ref{crankgen}, we have the following  relation:
\begin{equation}
\left<J_z\right>_q^{(2)}=C(z;q)
\end{equation}
We see $J_z(\lambda)$ and $z^{\operatorname{crk}(\lambda)}$ are related through a family of $q$-brackets; then using Corollaries \ref{class1}, \ref{class2}, and \ref{class3}, we can write $J_z(\lambda)$ explicitly. Noting that $J_z(\lambda)=(1-z)j_z(\lambda)$, we arrive at the formula we seek.
  
\begin{corollary}\label{jtpcoeff}
The partition-indexed coefficients of $j(z;q)^{-1}$ are
\[
j_z(\lambda)=(1-z)^{-1}\sum_{\delta|\lambda}\sum_{\varepsilon|\delta}z^{\operatorname{crk}(\varepsilon)}
\]
for $z\neq1$. In terms of the coefficients $M_z(*)$ given by Equation \ref{M}:
\[
j_z(\lambda)=(1-z)^{-1}\sum_{\delta|\lambda}\sum_{\substack{\varepsilon|\delta \\  \varepsilon\in\mathcal P_{=}}}\sum_{d|\ell(\varepsilon)}M_z(d)\mu\left(\frac{\ell(\varepsilon)}{d}\right)
\]
\end{corollary}

\begin{remark}
By Corollary \ref{two-way}, we can also write $M_z(n)$ in terms of the coefficients of $j(z;q)^{-1}$.
\end{remark}

We examine the function $j(z;q)^{-1}$ from a somewhat different perspective in another study \cite{SchneiderJTP}, connecting the $q$-bracket operator to the mock theta functions of Ramanujan.

\  
\ 
\  
\  
\  
\  
\  

\section*{Appendix:  \texorpdfstring{$q$}{Lg}-bracket arithmetic}
The $q$-bracket operator is reasonably well-behaved as an algebraic object; here we give a few formulas that may be useful for computation. 

Take $f,g,h:\mathcal P\to\C$. From Definition \ref{qbracket} we have $q$-bracket addition
\[
\left<f\right>_q+\left<g\right>_q=\left<f+g\right>_q
,\]
which is commutative, of course, and also associative:
\[
\left<f+g\right>_q+\left<h\right>_q=\left<f\right>_q+\left<g+h\right>_q
\]
We have for a constant $c\in\C$ that $c\left<f\right>_q=\left<cf\right>_q$; other basic arithmetic relations such as $\left<0\right>_q=0$ and $\left<f\right>_q+\left<0\right>_q=\left<f\right>_q$ follow easily as well.

Now, let us define $\tilde{f}:\Z_{\geq 0}\to\C$ by
\[
\tilde{f}(n):=\sum_{\lambda\vdash n}f(\lambda)
.\]
We define a convolution ``$*$'' of two such functions $\tilde{f},\tilde{g}$ by
\begin{equation}\label{convolution}
(\tilde{f}*\tilde{g})(\lambda):=\frac{1}{p(|\lambda|)}\sum_{k=0}^{|\lambda|}\tilde{f}(k)\tilde{g}(|\lambda|-k)
.
\end{equation}
Note that, by symmetry, $\tilde{f}*\tilde{g}=\tilde{g}*\tilde{f}$. 

Let us also define a multiplication ``$\star$'' between $q$-brackets by
\begin{equation}\label{times}
\left<f\right>_q\star\left<g\right>_q:=\frac{\left<f\right>_q\left<g\right>_q}{(q;q)_{\infty}}
,
\end{equation}
where the product and quotient on the right are taken in $\C[[q]]$. It follows from \ref{convolution} and \ref{times} above that
\[
\left<f\right>_q\star\left<g\right>_q=\left<\tilde{f}*\tilde{g}\right>_q 
.\]
From here it is easy to establish a $q$-bracket arithmetic yielding a commutative ring structure, with familiar-looking relations such as
\[\left<f\right>_q\star\left<\tilde{g}*\tilde{h}\right>_q=\left<\tilde{f}*\tilde{g}\right>_q\star\left<h\right>_q
,\]
\[
\left<f\right>_q\star\left<g+h\right>_q=\left<\tilde{f}*\tilde{g}\right>_q+\left<\tilde{f}*\tilde{h}\right>_q,
\]
and so on. 

It is trivial to see that $\left<1\right>_q=1$; however, $\left<1\right>_q\star\left<f\right>_q=\frac{\left<f\right>_q}{(q;q)_{\infty}}\neq\left<f\right>_q$, so $\left<1\right>_q$ is not the multiplicative identity in this arithmetic. In fact, as we note in Section \ref{section1}, the $q$-bracket of Dyson's rank function ``$\operatorname{rk}$'' is equal to $(q;q)_{\infty}
$. Then by Equation \ref{times}, $\left<\operatorname{rk}\right>_q$ may serve as multiplicative identity in the $q$-bracket arithmetic above (but is not unique in this respect, by the comments at the beginning of Section \ref{Applications of the $q$-antibracket}).

\ 
\ 
\ 

\section*{Acknowledgments}
The author is thankful to my colleagues Olivia Beckwith and Marie Jameson, for helping to clarify my arguments as I worked through these ideas, and to my Ph.D. advisor Ken Ono, for interesting discussions related to partition theory and the $q$-bracket operator. I am also very grateful to Krishnaswami Alladi, George Andrews, David Borthwick, John Duncan, Joel Riggs, Larry Rolen, and Andrew Sills, for stimulating conversations that informed this study. Furthermore, I wish to thank the anonymous referee for offering insightful comments and suggestions, and Tanay Wakhare for noticing an error in an earlier draft.  
\ 
\ 
\ 
\ 
\ 
\

\end{document}